%% file: hmgns_apprx.tex
\documentclass[11pt]{article}
\usepackage[dvipsnames]{xcolor}

\usepackage[utf8]{inputenc} 
\DeclareUnicodeCharacter{1EC7}{\d{\^{e}}}
\usepackage[T1]{fontenc} 

\usepackage[backend=biber, style=alphabetic, natbib=true]{biblatex} 
\usepackage[autostyle=true]{csquotes} 

\usepackage{mathpazo} 
\usepackage{amssymb}
\usepackage{amsthm}
\usepackage{mathtools}
\usepackage{xifthen}
\usepackage[shortcuts]{extdash} 
\usepackage{hyperref}
\usepackage{enumitem}
\usepackage{float}

\newtheorem{theorem}{Theorem}[section]
\newtheorem{proposition}[theorem]{Proposition}

\theoremstyle{definition}
\newtheorem{definition}[theorem]{Definition}
\newtheorem{remark}[theorem]{Remark}
\newtheorem{example}[theorem]{Example}

\usepackage[algoruled, linesnumbered]{algorithm2e}
\SetAlCapSty{textsc}
\SetKwSty{texttt}
\DontPrintSemicolon
\SetKw{True}{true}
\SetKw{False}{false}
\SetKw{Break}{break}
\SetKw{Unb}{unbounded}
\newcommand{\Assign}{\leftarrow}

\DeclareMathAlphabet{\mathbfcal}{OMS}{cmsy}{b}{n} 

\NewDocumentCommand{\set}{mg}{\left\lbrace {#1} \IfValueT{#2}{\,\middle|\, {#2}} \right\rbrace}
\AtBeginDocument{\renewcommand{\matrix}[2][]{\ifthenelse{\isempty{#1}}{\begin{pmatrix}
        #2 \end{pmatrix}}{\left( \begin{array}{#1} #2 \end{array} \right)}}}

\newcommand{\R}{\mathbb{R}}
\newcommand{\N}{\mathbb{N}}
\renewcommand{\S}{\mathcal{S}}
\newcommand{\A}{\mathcal{A}}
\newcommand{\B}{\mathcal{B}}
\newcommand{\T}{^{\mathsf{T}}}
\newcommand{\norm}[1]{\left\lVert {#1} \right\rVert}
\newcommand{\dist}[2]{d\left({#1},{#2}\right)}
\newcommand{\proj}[2][]{\ifthenelse{\isempty{#1}}{\pi\left({#2}\right)}{\pi_{#1}\left({#2}\right)}}
\newcommand{\hd}[2]{d_{\mathsf{H}}\left({#1},{#2}\right)} 
\newcommand{\thd}[2]{d_{\mathsf{tH}}\left({#1},{#2}\right)} 
\newcommand{\homog}[1]{\homogop {#1}} 
\newcommand{\rc}[1]{0^{\infty}{#1}}
\newcommand{\polar}[1]{{#1}^{\circ}}
\newcommand{\polarcone}[1]{{#1}^*}
\newcommand{\ed}{(\varepsilon,\delta)}
\newcommand{\decoRule}{\rule{.8\textwidth}{.4pt}}

\renewcommand{\geq}{\geqslant}
\renewcommand{\leq}{\leqslant}

\newcommand{\mgeq}{\succcurlyeq}

\newcommand{\mg}{\succ}

\DeclareRobustCommand{\paramone}{}
\DeclareRobustCommand{\paramtwo}{}
\DeclareRobustCommand{\paramthree}{}
\NewDocumentCommand{\varref}{mggg}{{\IfValueT{#2}{\renewcommand{\paramone}{#2}}\IfValueT{#3}{\renewcommand{\paramtwo}{#3}}\IfValueT{#4}{\renewcommand{\paramthree}{#4}}\ref{#1}}} 

\def\fchar#1#2@{#1} 
\def\estr#1#2@{#2} 
\NewDocumentEnvironment{optimization}{O{max} m o o s}
{
  \begin{equation*}
    \IfNoValueF{#3}{\ifthenelse{\equal{\fchar #3@}{*}}{\tag*{\estr #3@}}{\tag{#3}}} 
    \IfNoValueF{#4}{\label{#4}}
    \begin{aligned}
      \IfNoValueTF{#1}{\max\IfBooleanT{#5}{&}}{\ifthenelse{\equal{\detokenize{#1}}{\detokenize{min}}}{\min\IfBooleanT{#5}{&}}{\text{#1}\IfBooleanT{#5}{&}}}
      \; #2
      \IfBooleanTF{#5}{\\ \text{s.t.}&\;\begin{aligned}[t]}{\quad\text{s.t.}\quad}
}
{
  \IfBooleanT{#5}{\end{aligned}}
  \end{aligned}
  \end{equation*}
  \ignorespacesafterend
}

\makeatletter
\newcommand{\bcdot}{\mathpalette\scalebullet@{.5}}
\newcommand{\scalebullet@}[2]{\mathbin{\vcenter{\hbox{\scalebox{#2}{$\m@th#1\bullet$}}}}}
\makeatother

\makeatletter
\renewcommand*\env@matrix[1][*\c@MaxMatrixCols c]{%
  \hskip -\arraycolsep
  \let\@ifnextchar\new@ifnextchar
  \array{#1}}
\makeatother

\DeclareMathOperator{\cl}{cl}
\DeclareMathOperator{\interior}{int}
\DeclareMathOperator{\conv}{conv}
\DeclareMathOperator{\cone}{cone}
\DeclareMathOperator{\ext}{ext}
\DeclareMathOperator{\vertices}{vert}
\DeclareMathOperator{\homogop}{homog}
\DeclareMathOperator{\tr}{tr}

\newcommand{\figpath}{figures/}

\begin{document}

\title{Polyhedral approximation of spectrahedral shadows via
  homogenization}

\author{Daniel Dörfler \\
  Friedrich Schiller University Jena, Germany
  \and Andreas Löhne \\
  Friedrich Schiller University Jena, Germany}

\date{\today}

\maketitle

\begin{abstract}
  This article is concerned with the problem of approximating a not
  necessarily bounded spectrahedral shadow, a certain convex set, by
  polyhedra. By identifying the set with its homogenization the problem
  is reduced to the approximation of a closed convex cone. We introduce
  the notion of homogeneous $\delta$-approximation of a convex set and
  show that it defines a meaningful concept in the sense that
  approximations converge to the original set if the approximation error
  $\delta$ diminishes. Moreover, we show that a homogeneous
  $\delta$-approximation of the polar of a convex set is immediately
  available from an approximation of the set itself under mild
  conditions. Finally, we present an algorithm for the computation of
  homogeneous $\delta$-approximations of spectrahedral shadows and
  demonstrate it on examples.
\end{abstract}

\input{\secpath sec1}

\input{\secpath sec2}

\input{\secpath sec3}

\input{\secpath sec4}

\input{\secpath sec5}

{\noindent\footnotesize\textbf{Acknowledgements}\hspace{1em} The authors thank the two
  anonymous referees for their insightful comments which have enhanced
  the quality of the article.}

{\noindent\footnotesize\textbf{Data Availability Statement}\hspace{1em} Data sharing is
  not applicable to this article as no datasets were generated or
  analysed during the current study.}

\sloppy
\printbibliography[heading=bibintoc]
\fussy

\end{document}

%% file: sections/sec1.tex
\section{Introduction}\label{sec1}

Polyhedral approximation of convex sets is a profound method in convex
analysis with applications ranging through various branches of
mathematics. In the 1903 article \cite{Min03} Minkowski argues that
every compact convex set in three-dimensional Euclidean space can be
approximated arbitrarily well by \emph{polyhedra}, i.e. intersections
of finitely many closed halfspaces. Although it is not explicitly
stated in the article, his line of argumentation also applies to
general $n$-dimensional Euclidean spaces. This is later mentioned by
Bonnesen and Fenchel in \cite{BF34}.\\
Since then, polyhedral approximation has been applied in mathematical
programming, in particular in multiple objective optimization, e.g.
\cite{LRU14, DLSW22, WURKH23}, in approximation methods for general
convex optimization problems \cite{Vei67, BY11} and for mixed-integer
convex programs \cite{WP95, KLW16}. Other areas of application include
machine learning \cite{GLLS19, YBR08} and large deviation probability
theory \cite{NR95}.\\
In \cite{Kam92} Kamenev introduces a class of approximation algorithms
that compute outer and inner polyhedral approximations of a compact
convex set called cutting and augmenting schemes, respectively, by
iteratively improving upon a current approximation. The prevalent
method of quantifying the approximation error in these algorithms is
via the Hausdorff distance. Kamenev proves that under certain
conditions the sequence of polyhedral approximations generated by such
a scheme converges to the original set in the Hausdorff distance and
derives explicit bounds for the Hausdorff distance between an
approximation and the set to be approximated. These general schemes
subsume many of the polyhedral approximation algorithms in the
literature.\\
The literature concerned with polyhedral approximation of not
necessarily bounded sets is scarce. The reason appears to be the fact
that the Hausdorff distance only defines a metric on the family of
nonempty compact subsets of $\R^n$ but not on the larger family of
nonempty closed subsets. In fact, restrictive conditions need to be
satisfied for polyhedral approximation of unbounded sets in the
Hausdorff distance to be possible. For instance, the recession cones
of the approximation and the original set need to coincide, which
limits approximation to sets with polyhedral recession cones. A
characterization of those sets that are approximable by polyhedra in
the Hausdorff distance is due to Ney and Robinson \cite{NR95}. Similar
results are obtained in \cite{Ulu18} in the context of convex multiple
objective optimization. Evidently, the Hausdorff distance is no
suitable measure of approximation quality in the unbounded case. One
of the authors recently introduced the notion of
$\ed$\-/approximation, a concept of polyhedral outer approximation of
closed convex sets containing no lines, see \cite{Doe22}. One feature
is that these approximations bound the so-called truncated Hausdorff
distance, a metric on the family of closed cones in $\R^n$, between
the recession cones, but allow them to differ. Here, we present
another notion of polyhedral approximation that is not limited to
outer approximation, applicable to a larger class of sets and exhibits
elegant behavior under polarity not present for
$\ed$\-/approximations. These approximations, denoted homogeneous
$\delta$\-/approximations, are motivated by the fact that every closed
convex set in $\R^n$ can be identified with a closed convex cone in
$\R^{n+1}$, called its homogenization, and that it can also be
recovered therefrom. They are defined via the truncated Hausdorff
distance between the homogenizations of the approximation and the
original set.\\
This article is structured as follows. In the next section we
introduce necessary notation and basic definitions. Section \ref{sec3}
presents the notion of homogeneous $\delta$\-/approximation for the
polyhedral approximation of a convex set $C$. We show that homogeneous
$\delta$\-/approximations converge to $\cl C$ in the sense of
Painlevé-Kuratowski if the approximation error $\delta$ tends to
zero. Moreover, we prove that a polyhedron $P$ is a homogeneous
$\delta$\-/approximation of $C$ given that the Hausdorff distance
between $P$ and $C$ is finite and compute the approximation error
$\delta$ in terms of this distance. This shows that the concept of
homogeneous $\delta$\-/approximation is more general than polyhedral
approximation with respect to the Hausdorff distance. Finally, we
investigate the behaviour of homogeneous $\delta$\-/approximations
under polarity and show that a homogeneous $\delta$\-/approximation of
the polar set $\polar{C}$ can be obtained from a homogeneous
$\delta$\-/approximation of $C$ under mild assumptions. In the last
section an algorithm for the computation of outer and inner
homogeneous $\delta$\-/approximations for a special class of convex
sets called spectrahedral shadows is presented. These sets are the
images of spectrahedra, the feasible regions of semidefinite programs,
under linear transformations and are special in the sense that their
homogenizations can easily be described explicitly as spectrahedral
shadows at least up to closure. The algorithm is based on a recently
presented algorithm for the computation of polyhedral approximations
of the recession cone of a spectrahedral shadow, see \cite{DL22}. Two
examples illustrating the algorithm are presented.

%% file: sections/sec2.tex
\section{Preliminaries}\label{sec2}

Given a set $C \subseteq \R^n$ we denote by $\cl C$, $\interior C$ and
$\conv C$ the \emph{closure}, \emph{interior} and \emph{convex hull}
of $C$, respectively. The Euclidean unit ball of appropriate dimension
is written as $B$. A nonempty set $C$ is called a \emph{cone} if $\mu
x \in C$ for every $x \in C$ and $\mu \geq 0$, i.e. $C$ is closed
under nonnegative scalar multiplication. For convex $C$ its
\emph{recession cone} is $\set{d \in \R^n}{\forall x \in C, \forall
\mu \geq 0\colon x+\mu d \in C}$, denoted by $\rc{C}$. If $x \in \cl
C$, then $\rc{\cl C} = \set{d \in \R^n}{\forall \mu \geq 0\colon x+\mu
d \in \cl C}$, see \cite[Theorem 8.3]{Roc70}. The smallest, with
respect to inclusion, convex cone containing $C$ is called the
\emph{conical hull} of $C$ and written $\cone C$. We use the
convention $\cone\emptyset = \set{0}$. Furthermore, a point $x \in C$
is called an \emph{extreme point} of $C$ if $C \setminus \set{x}$ is
convex. The set of extreme points of $C$ is written $\ext C$. We
denote the hyperplane with normal vector $\omega \in \R^n$ and offset
$\gamma \in \R$ by $H(\omega, \gamma)$, i.e. $H(\omega,\gamma) =
\set{x \in \R^n}{\omega\T x = \gamma}$, and the lower closed halfspace
$\set{x \in \R^n}{\omega\T x \leq \gamma}$ by $H^-(\omega,\gamma)$. A
\emph{polyhedron} $P \subseteq \R^n$ is the intersection of finitely
many such halfspaces. Its extreme points are called \emph{vertices}
and are denoted by $\vertices P$. According to the well-known
Weyl-Minkowski theorem, see \cite[Theorem 14.3]{Gru07}, every
polyhedron can alternatively be expressed as the Minkowski sum of the
convex hull of finitely many points and the conical hull of finitely
many directions, i.e. $P =
\conv\set{v_1,\dots,v_m}+\cone\set{d_1,\dots,d_r}$ for suitable
$v_1,\dots,v_m \in \R^n$ and $d_1,\dots,d_r \in
\R^n\setminus\set{0}$. The \emph{polar} $\polar{C}$ and \emph{polar
cone} $\polarcone{C}$ of $C$ are defined as $\set{u \in \R^n}{\forall
x \in C\colon x\T u \leq \gamma}$ for $\gamma=1$ and $\gamma=0$,
respectively. It is easily verified that $\polar{\left(\cone
C\right)}=\polarcone{C}$.\\
We define $\dist{x}{C} = \inf\set{\norm{x-y}}{y \in C}$ as the
Euclidean distance between $x$ and $C$. If $C$ is closed and convex,
the infimum is uniquely attained, see e.g. \cite{HL01}. In this case,
the point at which the infimum is attained is called the
\emph{projection of $x$ onto $C$} and denoted by $\proj[C]{x}$. Now,
the \emph{Hausdorff distance} between sets $C_1, C_2 \subseteq \R^n$
is then expressed as
\[
  \hd{C_1}{C_2} = \max\set{\sup\set{\dist{x}{C_2}}{x \in C_1},
  \sup\set{\dist{x}{C_1}}{x \in C_2}}.
\]
The function $d_{\mathsf{H}}$ defines a metric on the class of
nonempty compact subsets of $\R^n$. For convex cones $K_1, K_2
\subseteq \R^n$ the \emph{truncated Hausdorff distance} between $K_1$
and $K_2$, denoted $\thd{K_1}{K_2}$, is defined as
\[
  \hd{K_1 \cap B}{K_2 \cap B}.
\]
That is, it is the usual Hausdorff distance between $K_1$ and $K_2$
restricted to $B$. Since cones contain the origin, $\thd{K_1}{K_2}
\leq 1$ always holds. The truncated Hausdorff distance provides a
metric on the class of closed convex cones in $\R^n$, see \cite{WW67}.

%% file: sections/sec3.tex
\section{Approximation concept for convex sets}\label{sec3}

In this section we present a concept of polyhedral approximation of a
convex set $C \subseteq \R^n$ that need not be bounded. The truncated
Hausdorff distance being a metric between closed convex cones is the
motivation for our approach. In order to work in a conic setting, we
assign to $C$ a convex cone from $\R^{n+1}$ called its homogenization.

\begin{definition}
  Let $C \subseteq \R^n$ be a convex set. The \emph{homogenization} or
  \emph{conification} of $C$, denoted $\homog{C}$, is the cone
  \[
    \cl\cone\left(C \times \set{1}\right).
  \]
\end{definition}

Every closed convex set $C \subseteq \R^n$ can be identified with its
homogenization because $C$ can be recovered from $\homog{C}$. In
particular,
\begin{equation}\label{Eq3.1}
  C = \set{x \in \R^n}{\matrix{x\\1} \in \homog{C}},
\end{equation}
i.e. intersecting $\homog{C}$ with the hyperplane $H(e_{n+1},1) = \{x
\in \R^{n+1} \mid x_{n+1}=1\}$ and projecting the result onto the
first $n$ variables yields the original set $C$. In the general case,
the right hand side set in Equation \eqref{Eq3.1} equals $\cl C$
because $\homog{C} = \homog{\cl C}$, cf. Proposition 3.8 below. Thus,
we can identify convex sets with their homogenizations up to
closures. Using this correspondence, one can work with convex cones
entirely. This is a standard tool in convex analysis to study problems
in a conic framework, see e.g. \cite{Roc70, RW98, Bri20}. Here, we use
it to reduce the problem of approximating convex sets to the problem
of approximating convex cones.

\begin{definition}
  For a convex set $C \subseteq \R^n$ and $\delta \geq 0$ a polyhedron
  $P \subseteq \R^n$ is called a \emph{homogeneous
    $\delta$\-/approximation} of $C$ if
  \[
    \thd{\homog{P}}{\homog{C}} \leq \delta.
  \]
\end{definition}

Clearly, if $P$ is a homogeneous $\delta$-approximation of $C$ it is
also one of $\cl C$ because $\homog{C}=\homog{\cl C}$. The notion of
homogeneous $\delta$\-/approximation can be understood as a relative
error measure between the involved sets. Consider a convex set $C
\subseteq \R^n$ and a polyhedron $P \subseteq \R^n$. If $x \in P$ with
$\dist{x}{C} \leq \varepsilon$, then the distance of the corresponding
direction of $\homog{P}$ to $\homog{C}$ obeys the relation
\[
  \dist{\frac{1}{\sqrt{x\T x + 1}}\matrix{x\\1}}{\homog{C}} \leq
  \frac{\varepsilon}{\sqrt{x\T x + 1}}.
\]
In particular, the distance depends inversely on $\norm{x}$. Hence, if
$P$ is a homogeneous $\delta$\-/approximation of $C$, then points of
one set that are far from the origin are allowed a larger distance to
the other set than points closer to the origin.\\
This observation is explained from the fact that the distinction
between points and directions of a closed convex set $C$ collapses
when transitioning to its homogenization. One has the identity
\begin{equation}\label{Eq3.2}
  \homog{C} = \cone\left(C \times \set{1}\right) \cup \left(\rc{C}
    \times \set{0}\right),
\end{equation}
see \cite[Theorem 8.2]{Roc70}, i.e. points as well as directions of
$C$ correspond to directions of $\homog{C}$. Thus, a convergent
sequence $\set{(x_k,\mu_k)\T}_{k \in \N}$ of directions of $\homog{C}$
does not necessarily yield a convergent sequence in $C$ itself. For
instance, if $\set{(x_k,\mu_k)\T}_{k \in \N}$ converges to a direction
$(\bar{x},0)\T$ and provided all $\mu_k$ are positive, then the
sequence $\set{\mu_k^{-1}x_k\T}_{k \in \N} \subseteq C$ is
unbounded. This means that closeness of homogenizations with respect
to the truncated Hausdorff distance does not imply closeness of the
sets themselves with respect to the Hausdorff distance. In fact, the
Hausdorff distance, which measures absolute distances, between a
homogeneous $\delta$\-/approximation $P$ of $C$ and $C$ may be
infinite for arbitrary $\delta \in (0,1]$.

\begin{example}\label{Ex3.3}
  Let $C = [0,\infty)$ and define $P_{\delta} =
  \left[0,\sqrt{1-\delta^2}/\delta\right]$ for $\delta \in (0,1]$.
  Then, $\rc{C} = C$ and $P_{\delta}$ is compact for every $\delta \in
  (0,1]$, but one has
  \[
    \begin{aligned}
      \thd{\homog{P_{\delta}}}{\homog{C}} &=
      \norm{\matrix{1\\0}-\proj[\homog{P_{\delta}}]{\matrix{1\\0}}}\\
      &=
      \norm{\matrix{1\\0}-\delta\sqrt{1-\delta^2}\matrix{\frac{\sqrt{1-\delta^2}}{\delta}\\1}}\\
      &= \delta,
    \end{aligned}
  \]
  The second equality is evident by taking into account that the
  projection of $(1,0)\T$ onto $\homog{P_{\delta}}$ is attained on the
  ray generated by $\left(\sqrt{1-\delta^2}/\delta,1\right)\T$. Hence,
  $P_{\delta}$ is a homogeneous $\delta$\-/approximation of $C$ with
  $\hd{P_{\delta}}{C} = \infty$.
\end{example}

A converse relation does hold.

\begin{proposition}\label{Prop3.4}
  Let $C \subseteq \R^n$ be a closed convex set and $P \subseteq \R^n$
  be a polyhedron. If $\hd{P}{C} \leq \varepsilon$, then $P$ is a
  homogeneous $\varepsilon$-approximation of $C$.
\end{proposition}
\begin{proof}
  Consider the alternative formula for the Hausdorff distance
  \[
    \hd{P}{C} = \inf\set{t > 0}{P \subseteq C + tB, C \subseteq P +
      tB},
  \]
  see \cite{RW98}. By assumption it holds $P \subseteq C+\varepsilon
  B$ and $C \subseteq P+\varepsilon B$. We only show the inclusion
  $\homog{P} \cap B \subseteq \homog{C} \cap B + \varepsilon B$. The other
  inclusion is derived identically by interchanging the roles of $P$
  and $C$. It holds
  \begin{equation}\label{Eq3.3}
    \cone\left(P \times \set{1}\right) \cap B \subseteq
    \cone\left((C+\varepsilon B) \times \set{1}\right) \cap B.
  \end{equation}
  Let $x$ be an element from the right hand side set. Then there exist
  $c \in C$, $\mu \in [0,1]$ and $d \in B$ such that
  $x=\mu(c+\varepsilon d,1)\T$. Define $\bar{x} =\proj[\homog
  C]{x}$. Then
  \[
    \norm{x-\bar{x}} = \dist{x}{\homog C} \leq \norm{x-\mu(c,1)\T}
    \leq \mu\varepsilon \leq \varepsilon.
  \]
  Moreover, $\norm{\bar{x}} \leq \norm{x} \leq 1$, because the
  projection mapping onto $\homog{C}$ is nonexpansive, see
  \cite[Proposition 3.1.3]{HL01}. Thus, $x = \bar{x}+(x-\bar{x}) \in
  \homog{C} \cap B+\varepsilon B$. Using \eqref{Eq3.3} we obtain
  the inclusion $\cone\left(P \times \set{1}\right) \cap B \subseteq
  \homog{C} \cap B + \varepsilon B$ and, because the right hand side set is
  closed, also $\homog{P} \cap B \subseteq \homog{C} \cap B +
  \varepsilon B$.
\end{proof}

The above proposition and Example \ref{Ex3.3} show that the concept of
homogeneous $\delta$\-/approximation is more general than
approximation with respect to the Hausdorff distance. Moreover,
homogeneous $\delta$-approximations provide a meaningful notion of
approximation in the sense that a sequence of homogeneous
$\delta$-approximations of $C$ converge to $C$ if $\delta$ tends to
$0$. The common notion of convergence in this scenario is that of
Painlevé-Kuratowski, a theory of convergence for sequences of sets
that applies to a broader class than Hausdorff convergence but
coincides with it in a compact setting. It is discussed in detail in
the book \cite{RW98} by Rockafellar and Wets and we adhere to the
notation therein. Let $\mathcal{N}_{\infty}$ and
$\mathcal{N}^\#_{\infty}$ denote the collections $\set{N \subseteq
  \N}{\N \setminus N \text{ is finite}}$ and $\set{N \subseteq \N}{N
  \text{ is infinite}}$ of subsets of $\N$, respectively. The set
$\mathcal{N}_{\infty}$ can be regarded as the set of subsequences of
$\N$ that contain all natural numbers larger than some $\bar{n} \in
\N$ and $\mathcal{N}^\#_{\infty}$ as the set of all subsequences of
$\N$. Obviously, $\mathcal{N}_{\infty} \subseteq
\mathcal{N}^\#_{\infty}$.

\begin{definition}
    For a sequence $\set{M_k}_{k \in \N}$ of subsets of $\R^n$ the
  \emph{outer limit} is the set
  \[
    \set{x \in \R^n}{\forall\;\varepsilon > 0,\;\exists\;N \in
      \mathcal{N}_{\infty}^\#,\;\forall\;k \in N \colon x \in M_k +
      \varepsilon B},
  \]
  denoted by $\limsup_{k \to \infty} M_k$, and the \emph{inner limit}
  $\liminf_{k \to \infty} M_k$ is the set
  \[
    \set{x \in \R^n}{\forall\;\varepsilon > 0,\;\exists\;N \in
      \mathcal{N}_{\infty},\;\forall\;k \in N \colon x \in M_k +
      \varepsilon B}.
  \]
\end{definition}

The outer limit of a sequence $\set{M_k}_{k \in \N}$ can be
paraphrased as the set of all points for which every open
neighbourhood intersects infinitely many elements $M_k$ of the
sequence. Similarly, the inner limit is the set for which the same is
true but for all elements $M_k$ beyond some $\bar{k} \in \N$. In
particular, both sets are closed.

\begin{definition}
  A sequence $\set{M_k}_{k \in \N}$ of subsets of $\R^n$ is said to
  \emph{converge to} a set $M \subseteq \R^n$ \emph{in the sense of
    Painlevé-Kuratowski} or \emph{PK-converge}, written $M_k \to M$,
  if
  \[
    M = \limsup_{k \to \infty} M_k = \liminf_{k \to \infty} M_k.
  \]
\end{definition}

Now, homogeneous $\delta$-approximations define a suitable notion of
approximation for closed convex sets because PK-convergence can be
characterized in terms of homogenizations and the truncated Hausdorff
distance.

\begin{proposition}[{cf. \cite[Corollary 4.47]{RW98}}]
  Let $\set{C_k}_{k \in \N}$ be a sequence of closed convex subsets of
  $\R^n$ and $C \subseteq \R^n$ be a closed convex set. Then the
  following are equivalent:
  \begin{enumerate}[label=(\roman*)]
  \item $C_k \to C$,
  \item $\thd{\homog{C_k}}{\homog{C}} \to 0$.
  \end{enumerate}
\end{proposition}

Hence, if $\set{P_k}_{k \in \N}$ is a sequence of homogeneous
$\delta_k$-approximations of $C$ and $\lim_{k \to \infty} \delta_k =
0$, then $\set{P_k}_{k \in \N}$ PK-converges to $C$.\\ We will now
emphasize another property of homogeneous $\delta$\-/approximations
that is not present for approximations in the Hausdorff
distance. Surprisingly, the polar set of a homogeneous
$\delta$-approximation of $C$ is a homogeneous $\delta$-approximation
of $\polar{C}$ provided the origin is contained in the sets. This is
different from the compact case. If $\hd{P}{C} \leq \delta$, then in
general it holds
\[
  \hd{\polar{P}}{\polar{C}} \leq \frac{\delta}{r_0(P)r_0(C)},
\]
where $r_0(P)$ and $r_0(C)$ are the radius of the largest Euclidean
ball centered at the origin that is contained in $P$ and $C$,
respectively, see \cite[Lemma 1]{Kam02}. In particular, the Hausdorff
distance between $\polar{P}$ and $\polar{C}$ does not only depend on
$\varepsilon$ but also on constants related to the geometry of $C$. In
order to establish the result, we need the a connection between the
polar cone of $\homog{C}$ and $\polar{C}$. The following result is
originally presented as \cite[Theorem 3.1]{BBW22} for closed convex
sets $C$. Here, we prove a slight generalization.

\begin{proposition}\label{Prop3.8}
  Let $C \subseteq \R^n$ be a convex set containing the origin. Then
  \[
    \polarcone{\homog{\left(C\right)}} = \cl\cone\left(\polar{C}
      \times \set{-1}\right) = -\homog{\left(-\polar{C}\right)}.
  \]
\end{proposition}
\begin{proof}
  The first equality is proved in \cite[Theorem 3.1]{BBW22} for $C$
  being closed and convex. However, it also holds for arbitrary convex
  sets $C$ because $\polar{C}=\polar{\left(\cl C\right)}$ and $\homog{C}
  = \homog{\cl C}$. The inclusion $\homog{C} \subseteq \homog{\cl C}$ is
  immediate from $C \subseteq \cl C$. For the other direction consider a
  point
  \[
    \mu\matrix{x\\1} \in \cone\left(\cl C \times \set{1}\right),
  \]
  that means we fix $\mu \geq 0$ and $x \in \cl C$. Thus, there exists
  a sequence $\set{x_k}_{k \in \N} \subseteq C$ converging to $x$. For
  every $k \in \N$ it holds $\mu\left(x_k\T,1\right)\T \in \cone\left(C
    \times \set{1}\right)$. Therefore, $\mu\left(x\T,1\right)\T \in
  \cl\cone\left(C \times \set{1}\right)$. This implies
  \[
    \cone\left(\cl C \times \set{1}\right) \subseteq \cl\cone\left(C
      \times \set{1}\right)
  \]
  and by taking closures
  \[
    \homog{\cl C} \subseteq \homog{C}.
  \]
  To show the second equation in the original statement we compute
  \[
    \begin{aligned}
      \cl\cone\left(\polar{C} \times \set{-1}\right) &= \cl
      \set{\mu\matrix{x\\-1} \in \R^{n+1}}{\mu \geq 0, x \in
        \polar{C}}\\
      &= \cl -\set{\mu\matrix{x\\1} \in \R^{n+1}}{\mu \geq 0, x \in
        -\polar{C}}\\
      &= - \cl\set{\mu\matrix{x\\1} \in \R^{n+1}}{\mu \geq 0, x \in
        -\polar{C}}\\
      &= -\homog{\left(-\polar{C}\right)}.
    \end{aligned}
  \]
  The penultimate line follows from \cite[Theorem 9.1]{Roc70}.
\end{proof}

\begin{figure}
  \centering
  \includegraphics[]{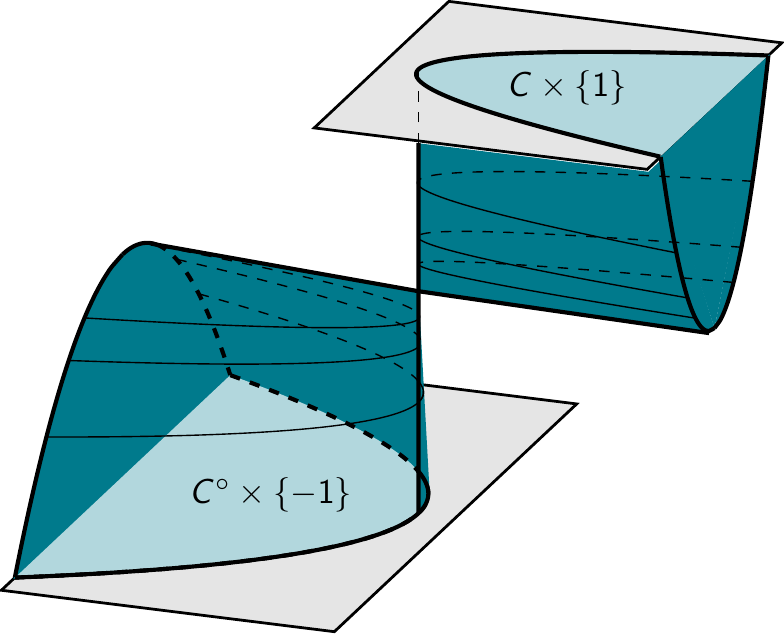}
  \decoRule
  \caption{Homogenization of the set
    $C=\set{x \in \R^2}{x_2 \geq x_1^2}$ and its polar cone. According
    to Proposition \ref{Prop3.8} the polar cone is the closure of the
    cone generated by $\polar{C} \times \set{-1}$.}
  \label{Fig3.1}
\end{figure}

Figure \ref{Fig3.1} illustrates Proposition \ref{Prop3.8} for the set
$C = \set{x \in \R^2}{x_2 \geq x_1^2}$.

\begin{theorem}\label{Th3.9}
  Let $C \subseteq \R^n$ be a convex set containing the origin. If $P$
  is a homogeneous $\delta$-approximation of $C$ containing the
  origin, then $\polar{P}$ is a homogeneous $\delta$-approximation
  of~$\polar{C}$.
\end{theorem}
\begin{proof}
  Let $M \colon \R^{n+1} \to \R^{n+1}$ denote the isometry defined by
  \[
    M\matrix{x\\\mu} = \matrix{x\\-\mu}.
  \]
  Now, we compute
  \[
    \begin{aligned}
      \thd{\homog{\left(\polar{P}\right)}}{\homog{\left(\polar{C}\right)}} &=
      \thd{M\left[\homog{\left(\polar{P}\right)}\right]}{M\left[\homog{\left(\polar{C}\right)}\right]}\\
      &= \thd{-\homog{\left(-\polar{P}\right)}}{-\homog{\left(-\polar{C}\right)}}\\
      &= \thd{\polarcone{\homog{\left(P\right)}}}{\polarcone{\homog{\left(C\right)}}}\\
      &= \thd{\homog{P}}{\homog{C}}\\
      &\leq \delta.
    \end{aligned}
  \]
  The first equation is true by \cite[Proposition 2.1]{IS10}. The second
  and third follow from Proposition \ref{Prop3.8} using the fact that
  $0 \in P \cap C$. Theorem 1 of \cite{WW67} yields the equality in line
  4. Finally, the inequality holds because $P$ is a homogeneous
  $\delta$-approximation of $C$. Thus, $\polar{P}$ is a homogeneous
  $\delta$-approximation of $\polar{C}$.
\end{proof}

%% file: sections/sec4.tex
\section{Computing homogeneous
\texorpdfstring{$\delta$}{d}-approximations of spectrahedral
shadows}\label{sec4}

In this section we present an algorithm for the computation of a
homogeneous $\delta$-approximation of a convex sets $C$. In order to
do this, we need to be able to explicitly describe the set
$\homog{C}$. Therefore, we restrict ourselves to a certain class of
convex sets for which this is possible called spectrahedral
shadows. Let $\S^{\ell}$ denote the linear space of real symmetric
$\ell \times \ell$-matrices. If $X \in \S^{\ell}$ is positive
semidefinite and positive definite, it is written as $X \mgeq 0$ and
$X \mg 0$, respectively. Moreover, for $X,Y \in \S^{\ell}$ an inner
product is defined by the trace of the matrix product, i.e. $X \bcdot
Y = \tr\left(XY\right)$. Now, a spectrahedral shadow $S \subseteq
\R^n$ can be represented as
\[
  S = \set{x \in \R^n}{\exists\;y \in \R^m\colon A_0 + \A(x) + \B(y)
    \mgeq 0}
\]
for $A_0 \in \S^{\ell}$ and linear functions $\A \colon \R^n \to
\S^{\ell}$, $\B \colon \R^m \to \S^{\ell}$. Every linear function $\A$
between $\R^n$ and $\S^{\ell}$ can be defined as $\A(x) = \sum_{i=1}^n
A_ix_i$ for suitable $A_1,\dots,A_n \in \S^{\ell}$, see
\cite{GM12}. Spectrahedral shadows are the images of spectrahedra, the
feasible regions of semidefinite programs, under orthogonal
projections. They are closed under many set operations,
e.g. intersections, Minkowski sums, linear transformations, taking
polar or conical hulls, see, for instance, \cite{GR95, HN12, Net11,
BPT13, Doe23}.\\
The algorithm we present is based on an algorithm for the polyhedral
approximation of recession cones of convex sets that was recently and
independently developed in \cite{Doe22} and \cite{WURKH23}. Given a
spectrahedral shadow $S$ and an error tolerance $\delta > 0$,
Algorithm 2 of \cite{DL22} computes polyhedral cones $K_{\mathsf{O}}$
and $K_{\mathsf{I}}$ such that $K_{\mathsf{I}} \subseteq \rc{\cl S}
\subseteq K_{\mathsf{O}}$, $\thd{K_{\mathsf{O}}}{\rc{\cl S}} \leq
\delta$ and $\thd{K_{\mathsf{I}}}{\rc{\cl S}} \leq \delta$. It is
shown that the algorithm works correctly and terminates after finitely
many steps under the following assumptions:
\begin{enumerate}[label=(A\arabic*)]
\item\label{Ass1} a point $p \in \R^n$ is knwon for which there
  exists $y_p \in \R^m$ such that $A_0+\A(x)+\B(y_p) \mg 0$,
\item\label{Ass2} a direction $\bar{d} \in \interior\rc{S} \cap B$
  is known.
\end{enumerate}
For ease of reference we recall the algorithm here as Algorithm
\ref{Alg1} below. Therefore, the following primal/dual pair of
semidefinite programs depending on parameters $p,d \in \R^n$ and a
spectrahedral shadow $S$ is needed.
\renewcommand{\paramone}{p}\renewcommand{\paramtwo}{d}\renewcommand{\paramthree}{S}
\begin{optimization}[max]{t}[*$\left(\text{P}\left(\protecting{\paramone},\protecting{\paramtwo},\protecting{\paramthree}\right)\right)$][P]
  A_0 + \A(x) + \B(y) &\mgeq 0\\
  x &= p+td
\end{optimization}
\begin{optimization}[min]{\left(A_0+\A(p)\right) \bcdot V}[*$\left(\text{D}\left(\protecting{\paramone},\protecting{\paramtwo},\protecting{\paramthree}\right)\right)$][D]
  -\A\T(V) &= w\\
  -\B\T(V) &= 0\\
  V &\mgeq 0\\
  d\T w &= 1
\end{optimization}
Here, $\A\T$ and $\B\T$ denote the adjoints of $\A$ and $\B$,
respectively. It holds $\A\T(V) = \left(A_1 \bcdot V,\dots,A_n \bcdot
V\right)\T$ and the corresponding for $\B\T$, see \cite[Lemma
4.5.3]{GM12}. Typically, we assume that $p \in S$ and $d$ is a
direction. Then problem \varref{P}{p}{d}{S} can be motivated as
follows. Starting at the point $p \in S$ the maximum distance is to be
determined that can be moved in direction $d$ without leaving the
set. If $S$ is compact, then the maximum distance will be attained in
the point where halfline $\set{p+td}{t \geq 0}$ intersects the
boundary of $S$. A solution to the dual program, if it exists, gives
rise to a supporting hyperplane of $\cl S$. The following proposition
is a slight variation of \cite[Proposition 3.11]{Doe23}.

\begin{proposition}\label{Prop4.1}
  Let $S \subseteq \R^n$ be a spectrahedral shadow and $p \in S$. Then
  the following hold:
  \begin{enumerate}[label=(\roman*)]
  \item\label{Prop4.1.i} if \varref{P}{p}{d}{S} is unbounded, then $d
    \in \rc{\cl S}$,
  \item\label{Prop4.1.ii} if \varref{P}{p}{d}{S} is strictly feasible
    with finite optimal value $t^*$, then an optimal solution $(V^*,w^*)$
    to \varref{D}{p}{d}{S} exists and the hyperplane
    $H(w^*,w^{*\mathsf{T}}p+t^*)$ supports $\cl S$ at $p+t^*d$.
  \end{enumerate}
\end{proposition}

\begin{algorithm}
  \caption{Approximation algorithm for recession cones of
    spectrahedral shadows}\label{Alg1}
  \KwIn{a spectrahedral shadow $S \subseteq \R^n$, point $p$
    satisfying \ref{Ass1}, direction $\bar{d}$ satisfying \ref{Ass2},
    error tolerance $\delta > 0$}
  \KwOut{outer and inner polyhedral approximation $K_{\mathsf{O}}$
    and $K_{\mathsf{I}}$ of $\rc{\cl S}$ with
    $\thd{K_{\mathsf{O}}}{K_\mathsf{I}} \leq \delta$}
  \eIf{\varref{P}{p}{-\bar{d}}{S} is \Unb}{
    $K_{\mathsf{O}} \Assign \R^n$\;
    $K_{\mathsf{I}} \Assign \R^n$\;
    \Return\;
  }{
    compute a solution $(V_{-\bar{d}},w_{-\bar{d}})$ to
    \varref{D}{p}{-\bar{d}}{S}\;
    $K_{\mathsf{O}} \Assign H^-(w_{-\bar{d}},0)$\;
    $K_{\mathsf{I}} \Assign \cone\set{\bar{d}}$\;
  }
  \Repeat{\upshape{stop}}{
    stop $\Assign$ \True\;
    \For{$r \in \vertices\set{x \in K_{\mathsf{O}}}{\norm{x}_{\infty}
        \leq 1} \setminus \set{0}$}{
      \For{$k=1$ \KwTo
        $\left\lceil\log_2\left(\frac{\norm{r-\bar{d}}}{\delta}\right)\right\rceil$}{
        $d \Assign \frac{2^k-1}{2^k}r + \frac{1}{2^k}\bar{d}$\;
        \eIf{\varref{P}{p}{d}{S} is \Unb}{
          $K_{\mathsf{I}} \Assign \cone\left(K_{\mathsf{I}} \cup \set{d}\right)$\;
        }{
          compute a solution $(V_d,w_d)$ to \varref{D}{p}{d}{S}\;
          $K_{\mathsf{O}} \Assign K_{\mathsf{O}} \cap H^-(w_d,0)$\;
          stop $\Assign$ \False\;
          \Break\;
        }
      }
    }
  }
\end{algorithm}

After finding initial outer and inner approximations $K_{\mathsf{O}}$
and $K_{\mathsf{I}}$ or determining that $S$ is the whole space in
lines 1--9, Algorithm \ref{Alg1} iteratively shrinks $K_{\mathsf{O}}$
or enlarges $K_{\mathsf{I}}$ until the truncated Hausdorff distance
between them is certifiably not larger than $\delta$. To this end, in
every iteration of the main loop in lines 10--25 a new search
direction $d$ is chosen in line 14 and the problem \varref{P}{p}{d}{S}
is investigated. If it is unbounded, then $d \in \rc{\cl S}$ according
to Proposition \ref{Prop4.1}\ref{Prop4.1.i} and $K_{\mathsf{I}}$ is
updated as $\cone\left(K_{\mathsf{I}} \cup \set{d}\right)$. Otherwise,
a solution $(V_d,w_d)$ to \varref{D}{p}{d}{S} exists and $H(w_d,w_d\T
p+t_d)$ is a supporting hyperplane of $\cl S$ for the optimal value
$t_d$ of \varref{P}{p}{d}{S} according to Proposition
\ref{Prop4.1}\ref{Prop4.1.ii}. Therefore, $H(w_d,0)$ supports $\rc{\cl
S}$ and the current outer approximation is updated as $K_{\mathsf{O}}
\cap H^-(w_d,0)$. After an update of either approximation occurs, a
new search direction is computed. For more details about the algorithm
we refer the reader to \cite{DL22}.\\
In order to apply Algorithm \ref{Alg1} to the problem of computing
homogeneous $\delta$-approximations of $S$, observe that every convex
cone is its own recession cone. This suggests we could apply Algorithm
\ref{Alg1} to $\homog{S}$ and obtain a homogeneous
$\delta$-approximation of $S$ by undoing the homogenization. This
approach presupposes that $S$ satisfies suitable assumptions such that
the requirements of Algorithm \ref{Alg1} are fulfilled for $\homog{S}$
and, more importantly, that $\homog{S}$ is a spectrahedral shadow for
which a representation is available or can be computed from a
representation of $S$. The following proposition shows that the latter
is possible.

\begin{proposition}\label{Prop4.2}
  Let $S = \set{x \in \R^n}{\exists\;y \in \R^m\colon A_0+\A(x)+\B(y)
    \mgeq 0}$. Then the homogenization $\homog{S}$ is the  closure of
  the spectrahedral shadow
  \[
    \set{\matrix{x\\x_{n+1}} \in \R^{n+1}}{\exists\;
      \begin{aligned}
        &y \in \R^m\\
        &t \in \R
      \end{aligned}\colon
      \begin{aligned}
        \A(x) + A_0x_{n+1} + \B(y) &\mgeq 0\\
        \matrix{x_{n+1} I & x\\x\T & t} &\mgeq 0
      \end{aligned}}.
  \]
\end{proposition}
\begin{proof}
  The Cartesian product $S \times \set{1}$ is the set $\{(x,x_{n+1})\T \in
  \R^{n+1}\mid\exists\;y\in \R^m\colon A_0+\A(x)+\B(y) \mgeq 0, x_{n+1} \geq
  1, -x_{n+1} \geq -1\}$. This set is a spectrahedral shadow because the
  three inequalities can be expressed as one matrix inequality using
  block diagonal matrices, i.e.
  \begin{align}\label{Eq4.1}
    S \times \set{1} &=
    \begin{multlined}[t]
      \left\lbrace\matrix{x\\x_{n+1}} \in \R^{n+1}\middle|\exists\;y \in
        \R^m\colon \matrix{A_0 & & \\ & -1 & \\ & & 1} + {} \right.\\
          \sum_{i=1}^n \matrix{A_i & & \\ & 0 & \\ & & 0}x_i +
          \left.\matrix{0 & & \\ & 1 & \\ & & -1}x_{n+1} + {} \right.\\
          \left.\sum_{i=1}^m \matrix{B_i & & \\ & 0 & \\ & & 0}y_i
            \mgeq \matrix{0 & & \\ & 0 & \\ & & 0}\right\rbrace.
    \end{multlined}
  \end{align}
  According to \cite[Proposition 4.3.1]{Net11} $\cone\left(S \times
    \set{1}\right)$ is the spectrahedral shadow
  \[
    \set{\matrix{x\\x_{n+1}} \in \R^{n+1}}{\exists\;
      \begin{aligned}
        &y \in \R^m\\
        &\mu \in \R\\
        &t \in \R
      \end{aligned}\colon
      \begin{aligned}
        \A(x) + A_0\mu + \B(y) &\mgeq 0\\
        x_{n+1} &\geq \mu\\
        -x_{n+1} &\geq -\mu\\
        \matrix{\mu & x_i\\x_i & t} &\mgeq 0 & i &= 1,\dots,n+1
      \end{aligned}}.
  \]
  Using the implicit equality $x_{n+1}=\mu$ we eliminate $\mu$ from
  the description of the set. The resulting spectrahedral shadow
  differs from the claimed set only through the occurence of the
  matrix inequalities
  \begin{equation}\label{Eq4.2}
    \matrix{x_{n+1} & x_i\\x_i & t} \mgeq 0
  \end{equation}
  for $i=1,\dots,n+1$. Now, there exists $t \in \R$ such that
  \eqref{Eq4.2} holds for $i=1,\dots,n+1$ if and only if the same is
  true but for $i=1,\dots,n$. Indeed, if $x_{n+1}=0$, then positive
  definiteness in \eqref{Eq4.2} implies $x=0$ and we can choose
  $t=0$ in both cases. Otherwise, if $x_{n+1} > 0$, then \eqref{Eq4.2}
  holds for some $i$, if and only if $t \geq x_{n+1}^{-1}x_i^2$
  according to the Schur complement, see \cite{Zha05}. Thus, we can
  set $t \geq x_{n+1}^{-1}\max\set{x_i^2}{i=1,\dots,n+1}$ in both
  cases. The case $x_{n+1} < 0$ does not occcur as it violates the
  positive semidefiniteness in \eqref{Eq4.2}. By a similar argument
  using Schur complement there exists $t \in \R$ such that
  \eqref{Eq4.2} holds for all $i=1,\dots,n$ if and only if there
  exists $t \in \R$ such that
  \[
    \matrix{x_{n+1}I & x\\x\T & t} \mgeq 0,
  \]
  where $I$ denotes the identity matrix of size $n$.
\end{proof}

Proposition \ref{Prop4.2} and its preceding paragraph motivate the
formulation of the following algorithm for the computation of
homogeneous $\delta$\-/approximations of $S$.

\begin{algorithm}[H]
  \caption{homogeneous $\delta$-approximation algorithm for
    spectrahedral shadows}\label{Alg2}
  \KwIn{a spectrahedral shadow $S \subseteq \R^n$, point $\bar{x} \in
    \interior S$ satisfying \ref{Ass1}, error tolerance $\delta > 0$}
  \KwOut{outer and inner homogeneous $\delta$-approximation
    $P_{\mathsf{O}}$ and $P_{\mathsf{I}}$ of $S$}
  $K \Assign \cone\left(S \times \set{1}\right)$\;
  compute outer and inner approximation $K_{\mathsf{O}}$ and
  $K_{\mathsf{I}}$ of $K$ using Algorithm \ref{Alg1} with input
  parameters $p=\left(\bar{x}\T,1\right)\T$, $\bar{d}=p/\lVert
  p\rVert$ and error tolerance $\delta$\;
  $P_{\mathsf{O}} \Assign \pi_{\R^n}\left[K_{\mathsf{O}} \cap
    H(e_{n+1},1)\right]$\;
  $P_{\mathsf{I}} \Assign \pi_{\R^n}\left[K_{\mathsf{I}} \cap
    H(e_{n+1},1)\right]$\;
\end{algorithm}

The function $\pi_{\R^n} \colon \R^{n+1} \to \R^n$ in lines 3 and 4
denotes the projection given by
\[
  \pi_{\R^n}\left(\matrix{x\\x_{n+1}}\right) = x.
\]

\begin{theorem}\label{Th4.3}
  Algorithm \ref{Alg2} works correctly, in particular it terminates
  with an outer homogeneous $\delta$-approximation $P_{\mathsf{O}}$ and
  an inner homogeneous $\delta$\-/approximation $P_{\mathsf{I}}$ of $S$.
\end{theorem}
\begin{proof}
  Since $\bar{x}$ satisfies Assumption \ref{Ass1} for $S$, the point $p
  = \left(\bar{x}\T,1\right)\T$ satisfies analogous conditions for $K$.
  This is seen from a direct calculation using the description of $K$
  derived in Proposition \ref{Prop4.2}. Similarly, $\bar{x} \in
  \interior S$ implies $\left(\bar{x}\T,1\right)\T \in \interior K$. In
  particular,
  \[
    \bar{d} = \frac{1}{\sqrt{\bar{x}\T \bar{x} + 1}}\matrix{\bar{x}\\1}
  \]
  satisfies Assumption \ref{Ass2} for $K$ because $K = \rc{K}$. From
  the correctness of Algorithm \ref{Alg1}, see \cite[Theorem 4.4]{DL22},
  it follows that line 2 of Algorithm \ref{Alg2} works correctly and
  returns polyhedral cones $K_{\mathsf{O}}$ and $K_{\mathsf{I}}$
  satisfying
  \begin{equation}\label{Eq4.3}
    K_{\mathsf{I}} \subseteq \rc{\cl{K}} = \cl K = \homog{S} \subseteq
    K_{\mathsf{O}}
  \end{equation}
  as well as $\thd{K_{\mathsf{O}}}{\homog{S}} \leq \delta$ and
  $\thd{K_{\mathsf{I}}}{\homog{S}} \leq \delta$.  For the sets
  $P_{\mathsf{O}}$ and $P_{\mathsf{I}}$ in lines 3 and 4 it holds
  \[
    \homog{P_{\mathsf{O}}} = K_{\mathsf{O}} \cap H^+(e_{n+1},0)
    \quad\text{and}\quad \homog{P_{\mathsf{I}}} = K_{\mathsf{I}}.
  \]
  The difference arises from the inclusions in \eqref{Eq4.3} and the
  fact that homogenizations are contained in the halfspace
  $H^+(e_{n+1},0)$ by definition. However,
  \[
    \thd{K_{\mathsf{O}} \cap H^+(e_{n+1},0)}{\homog{S}} \leq
    \thd{K_{\mathsf{O}}}{\homog{S}} \leq \delta
  \]
  holds due to the inclusion $\homog{S} \subseteq K_{\mathsf{O}} \cap
  H^+(e_{n+1},0)$. Therefore, $P_{\mathsf{O}}$ and $P_{\mathsf{I}}$
  are homogeneous $\delta$-approximations of $S$.\\
  The finiteness of Algorithm \ref{Alg2} is implied by the finiteness
  of Algorithm \ref{Alg1}, see \cite[Theorem 4.4]{DL22}.
\end{proof}

\begin{remark}
  The inner homogeneous $\delta$-approximation $P_{\mathsf{I}}$ of $S$
  computed by Algorithm \ref{Alg2} is always compact. Due to the
  containment of its homogenization $K_{\mathsf{I}}$ in the set
  $\cone\left(S \times \set{1}\right)$, one has $\mu > 0$ whenever
  $(x,\mu)\T \in K_{\mathsf{I}} \setminus \set{0}$. Thus, every
  direction of $K_{\mathsf{I}}$ corresponds to a point of
  $P_{\mathsf{I}}$, in particular
  \[
    \rc{P_{\mathsf{I}}} = \pi_{\R^n}\left[K_{\mathsf{I}} \cap
      H(e_{n+1},0)\right] = \set{0}
  \]
  according to Equation \eqref{Eq3.2}.
\end{remark}

We conclude this section with two examples.

\begin{figure}
  \centering
  \includegraphics[]{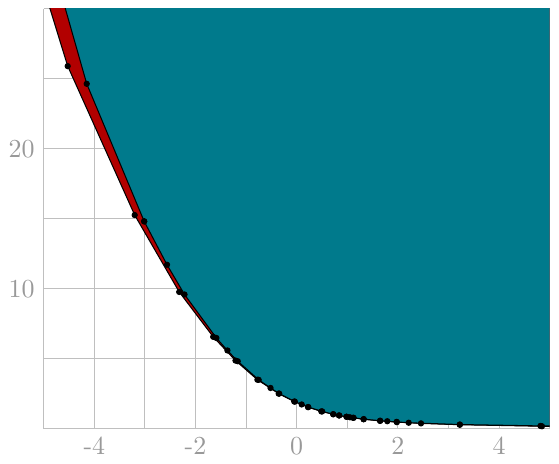}
  \includegraphics[]{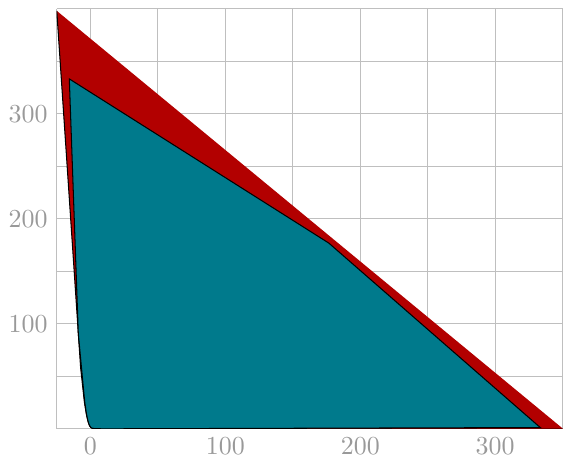}
  \decoRule
  \caption{Depicted are the outer and inner homogeneous
    $\delta$-approximation $P_{\mathsf{O}}$ and $P_{\mathsf{I}}$ of the
    set from Example \ref{Ex4.1} computed by Algorithm \ref{Alg2} with
    $\delta=0.01$ in red and blue, respectively. On the left is a section
    around the origin with vertices marked by black dots. Note that the
    distance between vertices of $P_{\mathsf{O}}$ and $P_{\mathsf{I}}$
    increases with their distance to the origin. The right figure shows
    the scale of $P_{\mathsf{I}}$, which is compact.}
  \label{Fig4.1}
\end{figure}

\begin{example}\label{Ex4.1}
  Algorithm \ref{Alg2} is applied to the spectrahedral shadow
  \[
    S = \set{x \in \R^2}{\exists\;y \in \R^2 \colon \matrix{x_1-y_1 &
        1\\1 & x_2-y_2} \mgeq 0, \matrix{y_2 & y_1\\y_1 & 1} \mgeq 0},
  \]
  which is the Minkowski sum of the epigraphs of the functions $x
  \mapsto x^{-1}$ restricted to $\R_{++}$ and $x \mapsto x^2$. As input
  parameters we set $\delta=0.01$ and $\bar{x} \approx
  (1.8846,1.8846)\T$. Figure \ref{Fig4.1} depicts $P_{\mathsf{O}}$ as
  well as the determined inner homogeneous $\delta$-approximation
  $P_{\mathsf{I}}$. In total, $534$ semidefinite programs are solved
  during the execution of Algorithm \ref{Alg2} for finding
  $P_{\mathsf{O}}$ and $P_{\mathsf{I}}$.
\end{example}

\begin{figure}[t]
  \centering
  \includegraphics[scale=.7]{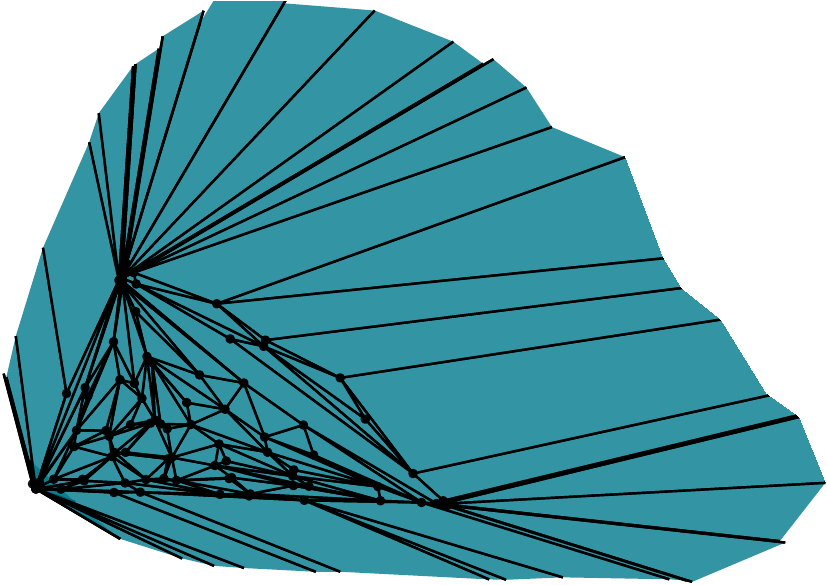}
  \hspace{.5cm}
  \includegraphics[scale=.5]{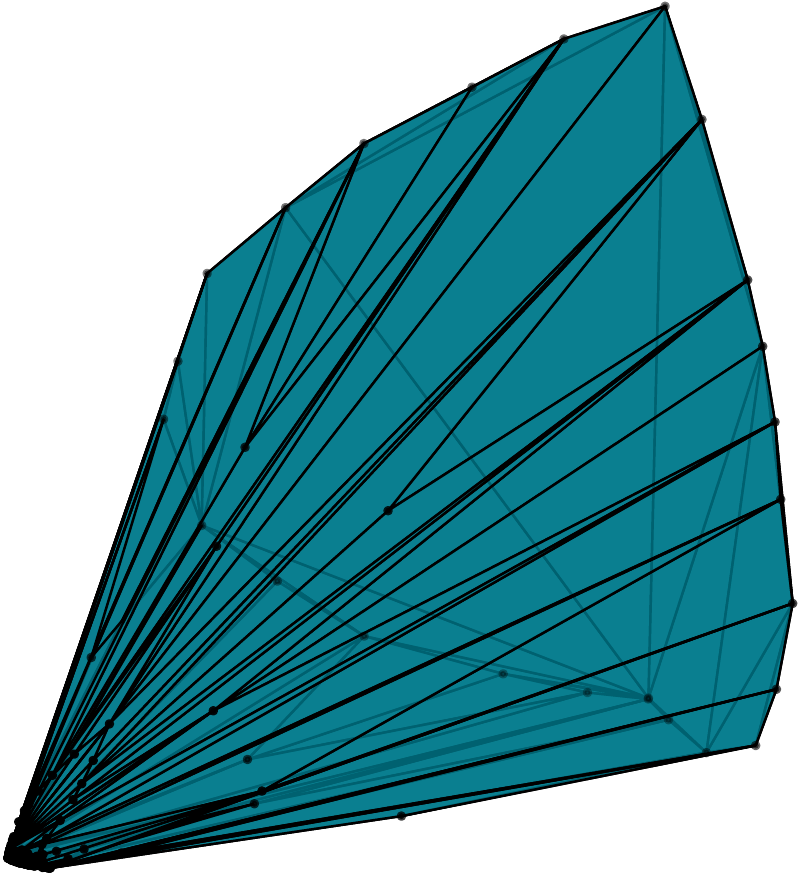}
  \includegraphics[scale=.65]{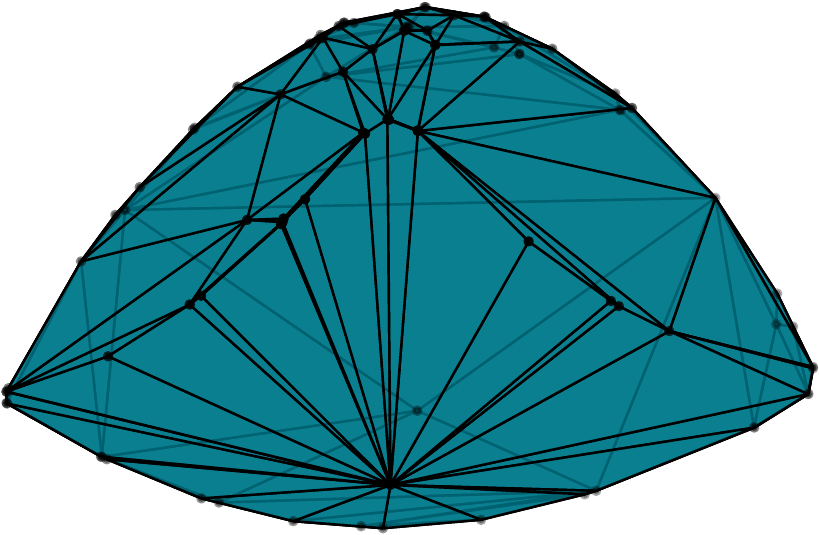}
  \hspace{.5cm}
  \includegraphics[scale=.61]{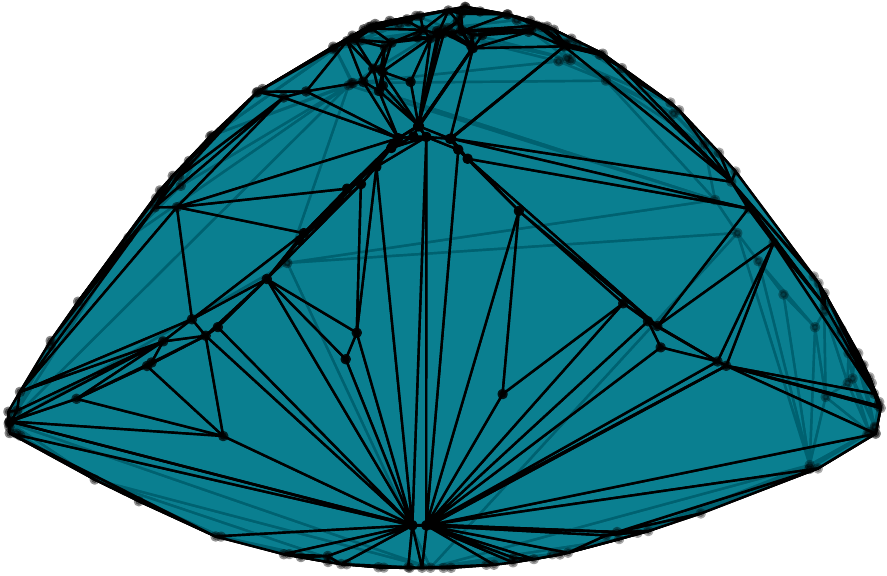}
  \decoRule
  \caption{The top row shows homogeneous $\delta$-approximations
    $P_{\mathsf{O}}$ and $P_{\mathsf{I}}$ of the set $S$ from Example
    \ref{Ex4.2} for $\delta=0.03$ restricted to balls of radius $7$ and
    $100$, respectively. Their polar sets, which are homogeneous
    $\delta$-approximations of $\polar{S}$, are depicted in the bottom
    row. Both are contained in a ball of radius $1.02$.}
  \label{Fig4.2}
\end{figure}

\begin{example}\label{Ex4.2}
  Consider the spectrahedral shadow $S$ consisting of points $x \in
  \R^3$ that admit the existence of $y^1, y^2 \in \R^3$ and $Z \in
  \S^3$, $Z \mgeq 0$, such that the system
  \[
    \begin{aligned}
      x-y^1-y^2 &= 0 & 2z_{13}+z_{22} &= 0\\
      z_{11}-1 &= 0 & z_{33}-1 &= 0\\
      \matrix{y^2_1+y^2_2 & 0 &
        2\left(y^2_1-y^2_2\right)\\0 & y^2_1+y^2_2 &
        2\sqrt{2}y^2_3\\2\left(y^2_1-y^2_2\right) & 2\sqrt{2}y^2_3 &
        y^2_1+y^2_2} &\mgeq 0 &
      \matrix{z_{11} & y^1_1 & y^1_2\\y^1_1 & z_{22} & y^1_3\\y^1_2 &
        y^1_3 & z_{33}} &\mgeq 0
    \end{aligned}
  \]
  is satisfied. Algorithm \ref{Alg2} with $\delta=0.03$ and
  $\bar{x}=0$ requires the solutions to $1989$ problems of type
  \varref{P}{p}{d}{K}, where $K=\cone\left(S \times
    \set{1}\right)$. The resulting outer approximation $P_{\mathsf{O}}$
  has $92$ vertices and $24$ extreme directions, while the compact inner
  approximation $P_{\mathsf{I}}$ is composed of $177$
  vertices. According to Theorem \ref{Th3.9} their polar sets
  $\polar{P}_{\mathsf{O}}$ and $\polar{P}_{\mathsf{I}}$ are an inner and
  outer homogeneous $\delta$-approximation of $\polar{S}$,
  respectively. They are both compact with $94$ and $232$ vertices
  each. Sets $P_{\mathsf{O}}$ and $P_{\mathsf{I}}$ as well as their
  polars are shown in Figure \ref{Fig4.2}.
\end{example}

%% file: sections/sec5.tex
\section{Conclusions}

We have introduced the notion of homogeneous $\delta$-approximation for
the polyhedral approximation of convex sets, which uses the common
concept of homogenization of a convex set. The advantage of this
approach compared to approximation with respect to the Hausdorff
distance is its compatibility with sets that are unbounded. Moreover,
we have shown that homogeneous $\delta$-approximations behave well under
polarity in the sense that the polar of a homogeneous
$\delta$-approximation of a convex set is a homogeneous
$\delta$-approximation of the polar set under mild assumptions. This
behaviour is not encountered when working with the Hausdorff
distance. Finally, we have presented an algorithm for the computation
of homogeneous $\delta$-approximations of spectrahedral shadows, a
subclass of convex sets. The algorithm is shown to be correct and
finite and its practicability is demonstrated on two examples.